\documentclass{article}
\usepackage{amsmath}[1996/11/01]
\usepackage{amssymb,amsthm,amsxtra}
\usepackage{epsfig}

\newlength{\mytopmargin}
\newlength{\myleftmargin}
\def\zz{\rlx\hbox{\small \sf Z\kern-.4em Z}}

\newtheorem{lemma}{Lemma}[section]

\newtheorem{prop}[lemma]{Proposition}
\newtheorem{cor}[lemma]{Corollary}

\begin{document}

\vspace{1cm}
\noindent
\begin{center}{   \large \bf On the
evaluation formula for Jack polynomials with \\
prescribed symmetry}
\end{center}
\vspace{5mm}

\noindent
\begin{center}
 P.J.~Forrester, D.S.~McAnally and Y.~Nikoyalevsky\\

\it Department of Mathematics and Statistics, \\
University of Melbourne, Victoria
3010, Australia
\end{center}
\vspace{.5cm}
\begin{quote}
The Jack polynomials with prescribed symmetry are obtained from
the nonsymmetric polynomials via the operations of symmetrization,
antisymmetrization and normalization. After dividing out the
corresponding antisymmetric polynomial of smallest degree, a
symmetric polynomial results. Of interest in applications is the value
of the latter polynomial when all the variables are set equal. Dunkl
has obtained this evaluation, making use of a certain skew symmetric
operator. We introduce a simpler operator for this purpose, thereby
obtaining a new derivation of the evaluation formula. 
An expansion formula of a certain product in terms of Jack polynomials
with prescribed symmetry implied by the evaluation formula is used
to derive a generalization of a constant term identity due to
Macdonald, Kadell and Kaneko.
Although
we don't give the details in this work, the operator introduced here
can be defined for any reduced crystallographic root system, and used to
provide an evaluation formula for the corresponding Heckman-Opdam
polynomials with prescribed symmetry.
\end{quote}

\vspace{.5cm}
\noindent
\section{Introduction}
\setcounter{equation}{0}
The type $A$ periodic
Calogero-Sutherland quantum many body system with exchange terms is
described by the Schr\"odinger operator
\begin{equation}\label{I.1}
H = - \sum_{j=1}^N {\partial^2 \over \partial x_j^2} +
\Big ( {\pi \over L} \Big )^2
\sum_{1 \le j < k \le N}
{(1/\alpha)((1/\alpha) - s_{jk}) \over
\sin^2(\pi (x_k - x_j)/L)}.
\end{equation}
In (\ref{I.1})
all particles are confined to the interval $[0,L]$ with
periodic boundary conditions, and $s_{jk}$ is the exchange operator
which permutes the variables $x_j$ and $x_k$. 
It is a standard result that
for
$0< x_1 < x_2 < \cdots < x_N <L$ the ground state for (\ref{I.1}) is
\begin{equation}\label{I.2}
\psi_0 := \prod_{1 \le j < k \le N}
\Big ( \sin ( \pi (x_k - x_j)/L) \Big )^{1/\alpha},
\end{equation}
and it is similarly a standard result that
with $z_j := e^{2 \pi i x_j / L}$ the excited states have the form
\begin{equation}\label{I.3}
\psi_0 \Big ( \prod_{j=1}^N z_j^l \Big ) f(z_1,\dots,z_N)
\end{equation}
where $l \in \mathbb Z_{\le 0}$ and $f$ is a homogeneous multivariable
polynomial. The polynomial $f$ is a type of Jack polynomial. 
In its most basic
form the Jack polynomial is nonsymmetric \cite{Op95}, 
however the transposition operator
$s_i := s_{i \, i+1}$ commutes with $H$ so we are free to restrict the
eigenfunctions to subspaces symmetric or antisymmetric with respect to
certain sets of coordinates. In physical terms this
corresponds to having a multicomponent system consisting of a mixture
of bosons and fermions. The polynomial part of the eigenfunction is then
referred to as a Jack polynomial with prescribed symmetry \cite{BF97b}.

The Jack polynomials with prescribed symmetry are to be denoted
$S_{\eta^*}^{(I,J)}(z)$, requiring the two labels $\eta^*$ and
$(I,J)$ for their unique specification.
The set  $I$ in the label $(I,J)$ determines the variables in which
$S_{\eta^*}^{(I,J)}$ is symmetric, while the set $J$ determines the
variables in which $S_{\eta^*}^{(I,J)}$ is antisymmetric.
Explicitly
\begin{equation}\label{2.I1}
s_i S_{\eta^*}^{(I,J)}(z) = S_{\eta^*}^{(I,J)}(z) \quad (i \in I)
\qquad
\qquad
s_j S_{\eta^*}^{(I,J)}(z) = -S_{\eta^*}^{(I,J)}(z) \quad (j \in J),
\end{equation}
and thus 
\begin{equation}\label{1.6a}
S_{\eta^*}^{(I,J)}(z) \propto {\rm Sym}_I {\rm Asym}_J \, E_\eta(z).
\end{equation}
For (\ref{2.I1}) to make sense we must have that
$I$ and $J$ are disjoint subsets of
$\{1,2,\dots,N-1\}$ such that 
\begin{equation}\label{es1}
i-1,i+1 \notin J \quad {\rm for} \quad i \notin I \quad {\rm and} \quad
j-1,j+1 \notin I \quad {\rm for} \quad
j \in J
\end{equation}
(thus representing the members of
$I$ by red dots and the members of $J$ by black dots on the lattice
$\{1,2,\dots,N-1\}$ there are no consecutive lattice points marked with
different coloured dots). The index $\eta^*$ is a composition such that
\begin{equation}\label{es}
\eta_i^* \ge \eta_{i+1}^* \quad \forall \quad i \in I, \qquad
\eta_j^* > \eta_{j+1}^* \quad \forall \quad j \in J.
\end{equation}

Let the set $J$ be decomposed as a union of sets of consecutive integers,
denoted $J_s$ say, and let $\tilde{J}_s := J_s \cup \{ \max(J_s) + 1 \}$
and
$\tilde{J} := \cup_s \tilde{J}_s$.
Then from (\ref{2.I1}) the polynomial $S_{\eta^*}^{(I,J)}$ can be
factorized in the form
\begin{equation}\label{2.I.1}
S_{\eta^*}^{(I,J)}(z) = \prod_s \Delta_{\tilde{J}_s}(z) \,
U_{\eta^*}^{(I,J)}(z), \quad
\Delta_X(z) := \prod_{i < i' \atop (i,i') \in X} (z_i - z_{i'}),
\end{equation}
where $U_{\eta^*}^{(I,J)}(z)$ is symmetric with respect to $s_i$ for 
$i \in I \cup J$.
Our interest is in the evaluation of $U_{\eta^*}^{(I,J)}(1^N)$, where
$$
1^N := (z_1,\dots, z_N) \Big |_{z_j=1 \, (j=1,\dots,N)}.
$$ 
In fact the value of the closely related quantity
$$
{{\rm Sym}_I {\rm Asym}_J \, E_\eta(z) \over
\prod_s \Delta_{\tilde{J}_s}(z)} \bigg |_{z = 1^N}
$$
is already known from the work of
Dunkl \cite{Du98}. Recalling (\ref{1.6a}) we see
the value of $U_{\eta^*}^{(I,J)}(z)$
follows once the proportionality constant in the former is
determined. In Proposition \ref{p2} below we will determine the
proportionality constant. However we will not then be done with the
problem. Rather we seek a self contained derivation, motivated
by our desire to obtain the evaluation formula for the analogue
of $U_{\eta^*}^{(I,J)}(1^N)$ in the case of Heckman-Opdam polynomials
with prescribed symmetry.
On this point, we 
recall (see e.g.~\cite{Op98,NW00}) that the Jack polynomials are the 
special case
of the  Heckman-Opdam polynomials corresponding to the type $A$ root
system; the latter can be constructed for all 
reduced crystallographic root
systems and are also the polynomial part of certain Schr\"odinger operators.
The derivation given here does indeed
permit a generalization to this more general
setting, although we will reserve the details until a later publication.

By way of further motivation, we point out that knowledge of the
evaluation of $U_{\eta^*}^{(I,J)}(1^N)$ has been of essential use in the
exact computation of retarded Green functions for spin generalizations of
the Hamiltonian (\ref{I.1}) \cite{KY98}. The reason is that the value of
$U_{\eta^*}^{(I,J)}(1^N)$ combined with the Cauchy product expansion
involving  the $S_{\eta^*}^{(I,J)}$ gives the expansion of the product
\begin{equation}\label{I.3.1}
\prod_{i \notin \tilde{J}} (1 - z_i)^{r-1}
\prod_s \prod_{j \in \tilde{J}_s} (1 - z_j)^{r - |J_s|} 
\end{equation}
in terms of $\{S_{\eta^*}^{(I,J)}\}$, which is one of the main technical
requirements in the computation of the retarded Green functions. The
expansion of (\ref{I.3.1}) and a corresponding constant term identity will
be discussed in Section 4. 

We will begin in Section 2 by revising essential properties of the Jack
polynomials with prescribed symmetry. In Section 3 we will introduce
an operator $ O_J$ which transforms
$S_{\eta^*}^{(I,J)}$ to be proportional to 
$S_{\eta^*}^{(I \cup J, \emptyset)}$.
We will show that the same operator acting on the RHS of 
the first equation in (\ref{2.I.1}),
followed by evaluation at $z=1^N$, acts only on $\prod_s \Delta_{\tilde{J}_s}
(z)$. By evaluating $ O_J( \prod_s \Delta_{\tilde{J}_s}(z))$ the value
of $U_{\eta^*}^{(I,J)}(1^N)$ follows. This is stated in Proposition \ref{p20}
below. In Section 4
we contrast our method with that of Dunkl. We also revise the expansion
of (\ref{I.3.1}) in terms of $\{S_{\eta^*}^{(I,J)}\}$, obtained from
knowledge of $U_{\eta^*}^{(I,J)}(1^N)$. This is then used to obtain a
certain constant term identity involving $S_{\eta^*}^{(I,J)}(z)$.

\section{Properties of the polynomials $S_{\eta^*}^{(I,J)}$}
\setcounter{equation}{0}
\subsection{Definitions and preliminaries}
Let $s_{jk}$ denote the permutation operator which when acting on
functions $f=f(z_1,\dots,z_N)$ interchanges $z_j$ and $z_k$. In the notation
of 
Knop and Sahi \cite{KS97} we introduce the type $A$
Cherednik operators by 
\begin{equation}\label{1.1}
\xi_i = \alpha z_i {\partial \over \partial z_i} +
\sum_{p < i} {z_i \over z_i - z_p} (1 - s_{ip}) +
\sum_{p>i} {z_p \over z_i - z_p} (1 - s_{ip}) +1 - i, \quad
i=1,\dots,N.
\end{equation}
The relation to the type $A$ root system becomes apparent by
comparing (\ref{1.1}) to the general Cherednik operator associated
with a root system,
\begin{equation}\label{1.2}
D_{\vec{\lambda}} :=
\sum_{j=1}^N \lambda_j z_j {\partial \over \partial z_j}
+ \sum_{\vec{\beta} \in R_+} k_{\vec{\beta}} { \langle \vec{\lambda} ,
\vec{\beta} \rangle \over z^{\vec{\beta}} - 1} (1 - s_{\vec{\beta}})
+ \langle \vec{\lambda}, \vec{\rho}_k \rangle,
\end{equation}
where $\vec{\lambda} = (\lambda_1,\dots,\lambda_N)$ denotes a real vector
with $N$-components,  $\langle \, , \,\rangle$ denotes the dot product  and
\begin{equation}\label{1.3}
z^{\vec{\beta}} := z_1^{\beta_1} z_2^{\beta_2} \cdots z_N^{\beta_N}.
\end{equation}
The quantities $R_+$, $ k_{\vec{\beta}}$, $s_{\vec{\beta}}$ and
$ \vec{\rho}_k$ depend on the root system $R$. In the type $A$ case, with
$\vec{e}_j$ denoting the elementary unit vector with component $j$
equal to $1$,
\begin{alignat}{2}
R_+ = \{ \vec{e}_j - \vec{e}_k | 1 \le j < k \le N \}, & \quad
k_{\vec{\beta}} = 1/\alpha \label{1.4} \\
s_{\vec{\beta}} = s_{jk} \: \: {\rm for} \: \: \vec{\beta} = 
\vec{e}_j - \vec{e}_k, & \quad \vec{\rho}_k = \sum_{j=1}^N
{1 \over \alpha} \Big ( {N+1 \over 2} - j \Big ) \vec{e}_j. \label{1.5}
\end{alignat}
Using these
explicit formulas in (\ref{1.1}), a short calculation shows that in the
type $A$ case
\begin{equation}\label{1.7}
D_{\vec{\lambda}} = {1 \over \alpha}
\sum_{j=1}^N \lambda_j \Big (\xi_j + {N - 1 \over 2} \Big ).
\end{equation}
In particular, with $\vec{\beta}$ as specified in (\ref{1.5}),
\begin{equation}\label{1.7a}
D_{\vec{\beta}} = {1 \over \alpha} (\xi_j - \xi_k).
\end{equation} 

A fundamental property of $\{ \xi_i \}$ is that they form a commuting
family of operators
which permit a set of
simultaneous polynomial eigenfunctions, labelled by a composition
$\eta := (\eta_1,\eta_2,\dots,\eta_N)$, $\eta_j \ge 0$, and 
homogeneous of degree $|\eta| := \sum_{j=1}^N \eta_j$. These 
are the nonsymmetric Jack polynomials $E_\eta(z;\alpha) =:
E_\eta(z)$, which are uniquely characterized as the
solution of the eigenvalue equation
\begin{equation}\label{1.7b}
\xi_i E_\eta = \bar{\eta}_i E_\eta, \qquad
\bar{\eta}_i := \alpha \eta_i - \#\{k < i | \eta_k \ge \eta_i \} -
\#\{k > i | \eta_k > \eta_i\}
\end{equation}
possessing a special triangularity structure when expanded in terms of
monomials. To specify this structure, let us denote by $\eta^+$ the
partition corresponding to the composition $\eta$. Let $\eta$ and
$\nu$, $\eta \ne \nu$ be compositions such that $|\eta| = |\nu|$, and
introduce the partial ordering $<$, known as the dominance ordering, by
the statement $\nu < \eta$ if
$\sum_{j=1}^p \nu_j < \sum_{j=1}^p \eta_j$ for each $p=1,\dots,N$.
A further partial ordering $\prec$ is defined on compositions by the
statement that $\nu \prec \eta$ if $\nu^+ < \eta^+$, or in the case
$\nu^+ = \eta^+$, if $\nu < \eta$. In terms of the monomials ordered by
the partial ordering $\prec$, the non-symmetric Jack polynomials have
the expansion
\begin{equation}\label{1.8}
E_\eta(z) = z^\eta + \sum_{\nu \prec \eta} c_{\eta \nu} z^\nu
\end{equation}
for some coefficients $c_{\eta \nu} = c_{\eta \nu}(\alpha)$.

Contained in the set of positive roots $R_+$ (recall (\ref{1.4})) are
the so called simple roots
\begin{equation}\label{3.1}
\vec{\alpha}_j := \vec{e}_j - \vec{e}_{j+1}, \qquad
j \in \{1,2,\dots,N-1\},
\end{equation}
which form a basis of $R_+$. We consider two subsets
$I,J \subseteq \{1,2,\dots,N-1\}$ such that $I \cap J = \emptyset$, and
write $W_{I \cup J} := \langle s_j | j \in I \cup J \rangle$ where
$s_j := s_{\vec{\alpha_j}} = s_{j \,j+1}$. We define the operation
${\cal O}_{I,J}$ on monomials by
$$
{\cal O}_{I,J}(z^\lambda) = \sum_{w \in W_{I \cup J}}
\det \nolimits_J (w) z^{w(\lambda)},
$$
where $w(\lambda) := (\lambda_{w^{-1}(1)},\dots, \lambda_{w^{-1}(N)})$,
and $\det_J(w) = 1$ if a decomposition of $w$ into a product of
transpositions has an even number of factors from
$\{s_j | j \in J\}$ and $\det_J(w) = -1$ if such a decomposition has an
odd number of entries from $\{s_j | j \in J\}$. This operator is extended
to general analytic functions by linearity.

Let us suppose that in addition to requiring $I \cap J = \emptyset$, we also
have the condition (\ref{es1}). In this circumstance
$W_{I \cup J} = W_I W_J = W_J W_I$ and furthermore
\begin{equation}\label{3.2}
{\cal O}_{I,J} = {\cal O}_I {\cal O}_J = {\cal O}_J {\cal O}_I
\quad {\rm with} \quad {\cal O}_I = {\rm Sym}_I = \prod_s {\rm Sym}_{I_s},
\: \:  {\cal O}_J = {\rm Asym}_J = \prod_s {\rm Asym}_{J_s} 
\end{equation}
where $I = \cup_s I_s$ and $J = \cup_s J_s$ (recall the beginning of the
paragraph below (\ref{2.I1})).

Our interest is in the polynomials 
${\cal O}_{I,J} E_\eta(z)$. From the definition of ${\cal O}_{I,J}$ we see
that 
\begin{equation}\label{4.2}
{\cal O}_{I,J} E_\eta(z) = 0, \quad
{\rm if} \: \: \eta_j = \eta_{j'} \: \: {\rm for \,\, any} \,\, j \ne j'
\in \tilde{J}_s, \,\, {\rm some} \,\, s,
\end{equation}
where $\tilde{J}_s$ is defined as in the beginning of the paragraph below
(\ref{es}).
On the other hand, if $\eta_j \ne  \eta_{j'}$ for all
$ j \ne j' \in \tilde{J}_s$, all $s$, then
\begin{equation}\label{4.3}
{\cal O}_{I,J} E_\eta(z) = a_\eta^{(I,J)} S_{\eta^*}^{(I,J)}(z), \qquad
S_{\eta^*}^{(I,J)}(z) = z^{\eta^*} + \sum_{\nu \prec \eta^*}
\tilde{c}_{\eta \nu} z^\nu
\end{equation}
for some $a_\eta^{(I,J)}$ and
$\tilde{c}_{\eta \nu} = \tilde{c}_{\eta \nu}(\alpha)$. In
(\ref{4.3}) $\eta^*$ is the unique element of $W_{I \cup J}(\eta)$ such that
(\ref{es}) holds. The polynomials
$S_{\eta^*}^{(I,J)}$, which were first introduced in \cite{BF97b}, and
subsequently studied in the works \cite{BDF00,Du98,KY98}, are referred to
as the Jack polynomials with prescribed symmetry. They have the symmetry
properties (\ref{2.I1}), 
which together with the structural formula in (\ref{4.3})
can be used to uniquely characterize the polynomials.

Our subsequent calculations require the explicit value of the
proportionality constant $a_\eta^{(I,J)}$ in (\ref{4.3}). One way
of obtaining this is to make use of the explicit expansion of
$S_{\eta^*}^{(I,J)}(z)$ in terms of $\{E_\nu\}$. Let us first revise
the derivation of this latter expansion \cite{BDF00}. Now, we know
\cite{Op95} that the action of $s_i$ on $E_\eta$ is given by
\begin{equation}\label{5.1}
s_i E_\eta(z) = \left\{ \begin{array}{ll}
{1 \over \bar{\delta}_{i,\eta}} E_\eta(z) + \Big ( 1 - {1 \over
\bar{\delta}_{i,\eta}^2} \Big ) E_{s_i \eta}(z), & \eta_i > \eta_{i+1} \\[.2cm]
E_\eta(z), & \eta_i = \eta_{i+1} \\[.2cm]
{1 \over \bar{\delta}_{i,\eta}} E_\eta(z) +
E_{s_i \eta}(z), & \eta_i < \eta_{i+1}
\end{array} \right.
\end{equation}
where $\bar{\delta}_{i,\eta} := \bar{\eta}_i - \bar{\eta}_{i+1}$, with
$\bar{\eta}_i$ specified by (\ref{1.7}). It follows from (\ref{5.1}),
(\ref{3.2}) and (\ref{4.3}) that
\begin{equation}\label{5.2}
S_{\eta^*}^{(I,J)}(z) = \sum_{\mu \in W_{I \cup J}(\eta^*)}
\hat{c}_{\eta^* \mu} E_\mu(z), \quad \hat{c}_{\eta^*\eta^*} = 1.
\end{equation}
Moreover, the coefficients $\hat{c}_{\eta^* \mu}$ in (\ref{5.2}) can
be computed explicitly in terms of the quantities
\begin{equation}\label{5.3}
d_\eta  := \prod_{(i,j) \in \eta} \Big ( \alpha(a(i,j) + 1) +
l(i,j) + 1 \Big ), \quad
d_\eta'  :=  \prod_{(i,j) \in \eta} \Big ( \alpha(a(i,j) + 1) +
l(i,j) \Big ),
\end{equation}
where
the notation $(i,j) \in \eta$ refers to the diagram of
the composition $\eta$, while
\begin{equation}\label{5.3a}
a(i,j) := \eta_i - j, \qquad
l(i,j) := \#\{k<i | j \le \eta_k + 1 \le \eta_i \} +
\# \{ k > i | j \le \eta_k \le \eta_i \}.
\end{equation}
Our derivation makes use of the fact that
the quantities $d_\eta$ and $d_{\eta'}$ have the properties 
\cite{Sa96}
\begin{equation}\label{6.1}
{d_{s_i \eta} \over d_\eta} = {\bar{\delta}_{i,\eta} + 1 \over
\bar{\delta}_{i,\eta}}, \quad 
{d_{s_i \eta}' \over d_\eta'} =
{\bar{\delta}_{i,\eta} \over \bar{\delta}_{i,\eta} - 1},
\qquad \eta_i > \eta_{i+1}
\end{equation}

\begin{prop}\label{p1}
Let $w \in W_{I \cup J}$ be decomposed as $w = w_I w_J$ where
$w_I \in W_I$, $w_J \in W_J$. Let $w \eta^* = \mu$ and
$w_I \eta^* = \mu_I$. Then the coefficents in (\ref{5.2}) are
specified by
\begin{equation}\label{6.2}
\hat{c}_{\eta^* \mu} = \det\nolimits_J(w) {d_{\eta^*}' d_\mu \over
d_{\mu_I}' d_{\mu_I}}.
\end{equation}
\end{prop}

\noindent Proof. \quad 
We write 
\begin{equation}\label{6.3}
\sum_{\mu \in W_{I \cup J}(\eta^*)}
\hat{c}_{\eta^* \mu} E_\mu(z) =
\sum_{\mu \in W_{I \cup J}(\eta^*) \atop \mu_i \le \mu_{i+1}}
\chi_{\mu_i \mu_{i+1}} \Big ( \hat{c}_{\eta^* \mu} E_\mu(z) +
 \hat{c}_{\eta^* \, s_i\mu} E_{s_i \mu}(z) \Big )
\end{equation}
where $\chi_{\mu_i \mu_{i+1}} = {1 \over 2}$ for
$\mu_i = \mu_{i+1}$, $\chi_{\mu_i \mu_{i+1}} = 1$ otherwise.
Applying the operator $s_i$ for $i \in I$ and $s_j$ for $j \in J$
to (\ref{6.3}) with $\hat{c}_{\eta^* \mu}$ given by (\ref{6.2}), we see
by making use of (\ref{6.1}) that the properties (\ref{2.I1}) hold.
Also, with
$w_I=w_J = {\rm Id}$ (the identity), and thus
$\mu = \mu_I = \mu_J = \eta^*$ and $\det\nolimits_J(w) = 1$,
we see that $\hat{c}_{\eta^* \eta^*} = 1$.
\hfill $\square$

\medskip
Proposition \ref{p1} can be used to determine $a_\eta^{(I,J)}$. To present
the result requires some notation. Let 
$$
{\cal M}_{I,\eta} = \#\{ \sigma' \in W_I |
\sigma'(\eta ) = \eta \}
$$
and write $\eta^{(\epsilon_I, \epsilon_J)}$, where
$\epsilon_I, \epsilon_J \in \{+,0,-\}$ to denote the element of
$W_{I \cup J}(\eta)$ with the property that $\eta^{(+,\cdot)}$
($\eta^{(\cdot,+)}$) has $\eta^{(+,\cdot)}_i \ge \eta^{(+,\cdot)}_{i+1}$
for all $i \in I$ 
($\eta^{(\cdot,+)}_j > \eta^{(\cdot,+)}_{j+1}$
for all $j \in J$), $\eta^{(-,\cdot)}$ ($\eta^{(\cdot,-)}$) has
$\eta^{(-,\cdot)}_i \le \eta^{(-,\cdot)}_{i+1}$ for all $i \in I$
($\eta^{(\cdot,-)}_j < \eta^{(\cdot,-)}_{j+1}$
for all $j \in J$), while $\eta^{(0,\cdot)}$ $(\eta^{(\cdot,0)})$ has
$\eta^{(0,\cdot)}_i = \eta_i$, $\eta^{(0,\cdot)}_{i+1} = \eta_{i+1}$ 
for all $i \in I$
($\eta^{(\cdot,0)}_j = \eta_j$, $\eta^{(\cdot,0)}_{j+1} = \eta_{j+1}$ 
for all $j \in J$).

\begin{prop}\label{p2}
With $a_\eta^{(I,J)}$ defined by (\ref{4.3}) we have
\begin{equation}\label{8.1}
a_\eta^{(I,J)} = \det \nolimits_J(w_J) {\cal M}_{I, \eta}
{d_\eta' d_{\eta^{(-,+)}}' d_{\eta^{(-,+)}} \over
 d_{\eta^{(0,+)}}' d_{\eta^{(0,+)}}  d_{\eta^{(-,-)}}' }
\end{equation}
where $w_J$ is such that $w_J \eta = \eta^{(0,+)}$.
\end{prop}

\noindent
Proof. \quad We model our proof on the derivation given in \cite{BF97b}
of the value of $a_\eta^{(I,J)}$ in the case that $I = \{1,\dots,N-1\}$
and $J = \emptyset$ (symmetrization in all variables). The first step is
to introduce the polynomial
\begin{equation}\label{t1}
G(x,y) = \sum_{\mu \in W_{I \cup J}(\eta^*)}
{d_\mu \over d_\mu'} E_\mu(x) E_\mu(y),
\end{equation}
which from (\ref{5.1}) is seen to have the property
$$
s_i^{(x)} G(x,y) = s_i^{(y)} G(x,y) \qquad {\rm for} \quad i \in I \cup J.
$$
In particular
\begin{equation}\label{t2}
{\cal O}^{(x)} G(x,y) = {\cal O}^{(y)} G(x,y).  
\end{equation}
Substituting (\ref{t1}) in (\ref{t2}) and recalling (\ref{4.3}) shows
$$
S_{\eta^*}^{(I,J)}(x) \sum_{\mu \in W_{I \cup J}(\eta^*)}
{d_\mu \over d_\mu'} a_\mu^{(I,J)} E_\mu(y) =
S_{\eta^*}^{(I,J)}(y) \sum_{\mu \in W_{I \cup J}(\eta^*)}
{d_\mu \over d_\mu'} a_\mu^{(I,J)}  E_\mu(x).
$$
It follows from this that
\begin{equation}\label{9.1}
S_{\eta^*}^{(I,J)}(x) = a_{\eta^*} \sum_{\mu \in W_{I \cup J}(\eta^*)}
{d_\mu \over d_\mu'} a_\mu^{(I,J)} E_\mu(x)
\end{equation}
for some constant $ a_{\eta^*}$. Comparing (\ref{9.1}) with (\ref{5.2})
and recalling the result of Proposition \ref{p1} we thus have
\begin{equation}\label{9.2}
a_{\eta^*} {d_\mu \over d_\mu'}  a_\mu^{(I,J)} =
\det\nolimits_J(w_J) {d_{\eta^*}' d_\mu \over
d_{\mu_I}' d_{\mu_I}}
\end{equation}
with $\mu_I = \mu^{(+,0)}$ and $w_J$ such that
$w_J \mu = \mu^{(0,+)}$. The identity
(\ref{9.2}) must hold for all $\mu \in W_{I \cup J}(\eta^*)$.
In particular it must hold for the smallest composition with respect 
to the partial ordering $\prec$, namely $\eta^{(-,-)}$. This
composition is special in that
$$
{\cal O}_{I,J} E_{\eta^{(-,-)}}(z) = \det\nolimits_J(w_{J'})
{\cal M}_{I,\eta} S_{\eta^*}^{(I,J)}(z)
$$
where $w_{J'}\eta^{(-,-)} = \eta^{(-,+)}$. Thus
$$
a_{\eta^{(-,-)}}^{(I,J)} = \det\nolimits_J(w_{J'}) {\cal M}_{I,\eta},
$$
which when substituted in (\ref{9.2}) with $\mu = \eta^{(-,-)}$ implies
\begin{equation}\label{9.3}
a_{\eta^*} = {1 \over {\cal M}_{I,\mu} }
{d_{\eta^*}' d_{\eta^{(-,-)}}' \over d_{\eta^{(-,+)}}' d_{\eta^{(-,+)}}}
\end{equation}
Substituting (\ref{9.3}) in (\ref{9.2}) with $\mu = \eta$ (of course
$\eta \in W_{I \cup J}(\eta^*)$) gives
(\ref{8.1}).
\hfill $\square$
 
In our derivation of the evaluation formula for the polynomial
$U_{\eta^*}^{(I,J)}(z)$ in
(\ref{2.I.1}), it is a corollary of Proposition \ref{p2} giving the
evaluation formula for $S_{\eta^*}^{(I,\emptyset)}$ which is of use.

\begin{cor}\label{p10}
For the Jack polynomial with prescribed symmetry,
$S_{\eta^*}^{(I,\emptyset)}$,
 constructed out of
symmetrization operations only, we have
the evaluation formula
\begin{equation}
S_{\eta^*}^{(I,\emptyset)}(1^N) = 
{| {\rm Sym} \, I | \over a_{\eta^*}^{(I,\emptyset)}} E_{\eta^*}(1^N) =
{| {\rm Sym} \, I | \over
{\cal M}_{I,\eta^*}} {e_{\eta^+} \over d_{\eta^{*(-,\cdot)}} }
\end{equation}
where $\eta^+$ denotes the unique partition which can be formed from
$\eta^*$ and
\begin{equation}\label{se}
e_{\eta^+} := \alpha^{|\eta|} [1 + N/\alpha]_{\eta^+}^{(\alpha)}, \qquad
[u]_{\eta^+}^{(\alpha)} := \prod_{j=1}^N
{\Gamma(u - {1 \over \alpha} (j-1) + \eta_j^+) \over
\Gamma(u - {1 \over \alpha} (j-1) )}.
\end{equation}
\end{cor}

\noindent
Proof. \quad The first equality follows immediately from (\ref{4.3}) with
$J = \emptyset$, while the second equality follows from
Proposition \ref{p2} with $J =  \emptyset$ and $\eta = \eta^*$, which gives
$$
 a_\eta^{(I,\emptyset)} = {\cal M}_{I,\eta} {d_{\eta^{(-,\cdot)}} \over
d_\eta},
$$
and the well known result \cite{Sa96}
\begin{equation}\label{2.28a}
E_\eta(1^N) = { e_{\eta^+} \over d_\eta }.
\end{equation}
\hfill $\square$

\section{Evaluation formula for $U_{\eta^*}^{(I,J)}$}
\setcounter{equation}{0}
\subsection{A special operator}
Of central importance to our eventual evaluation of 
$U_{\eta^*}^{(I,J)}$ is the operator
\begin{equation}\label{11.0}
O_J := \prod_{\vec{\beta} \in R_{J,+}} \Big ( D_{\vec{\beta}} + {1 \over
\alpha} \Big )
\end{equation}
where 
$$
R_{J,+} = \{ \vec{e}_j - \vec{e}_k | 1 \le j < k \le N, \, j,k \in \tilde{J}_s,
\, {\rm some} \ s \}
$$
and $ D_{\vec{\beta}}$ is specified in terms of the $\xi_i$ by (\ref{1.7})
(in the case $R_{J,+} = R_+$ this operator is introduced for general reduced
crystallographic root systems in \cite{Op98}).
Also, $R_{J,-}$ is defined by 
$$
R_{J,-} = \{ - \vec{e}_j + \vec{e}_k | 1 \le j < k \le N, \, 
j,k \in \tilde{J}_s, \, {\rm some} \ s \}.
$$
We seek an algebraic relation for $s_i O_J$. First,
one can check from the relations
$$
\xi_i s_i - s_i \xi_{i+1} = 1, \quad
\xi_{i+1} s_i - s_i \xi_{i} = -1, \quad
[\xi_i, s_j ] = 0 \: \: (j \ne i-1,i)
$$
(the subalgebra of the degenerate type $A$ Hecke algebra satisfied by
$\{\xi_i, s_j \}$)
and (\ref{1.7a}), or alternatively directly from (\ref{1.2}), that
for $\vec{\beta}_i = \vec{e}_i - \vec{e}_{i+1}$,
\begin{equation}\label{11.1}
s_i D_{\vec{\beta}} - D_{s_i(\vec{\beta})} s_i = {1 \over \alpha}
\langle \vec{\beta} , \vec{\beta}_i \rangle. 
\end{equation}
Now we rewrite
(\ref{11.0}) to read
\begin{equation}\label{11.2}
O_J = \prod_{\vec{\beta} \in R_{J,+} \atop
\langle \vec{\beta} , \vec{\beta}_i \rangle = 0}
\Big ( D_{\vec{\beta}} + {1 \over \alpha} \Big )
\prod_{\vec{\beta} \in R_{J,+} \atop
\langle \vec{\beta} , \vec{\beta}_i \rangle \ne 0, \,
 \vec{\beta} \ne  \vec{\beta}_i}
\Big ( D_{\vec{\beta}} + {1 \over \alpha} \Big )
\left \{ \begin{array}{ll} 1, &  \vec{\beta}_i \notin R_{J,+} \\
\Big ( D_{\vec{\beta}_i} + {1 \over \alpha} \Big ), &
 \vec{\beta}_i \in R_{J,+} \end{array} \right.
\end{equation}
Denoting the first product of operators in (\ref{11.2}) by
$O_J^{(1)}$, we see immediately from (\ref{11.1}) that
\begin{equation}\label{O1}
s_i O_J^{(1)} = O_J^{(1)} s_i.
\end{equation}
For the second product of operators, $O_J^{(2)}$ say, we note that if
$\vec{\beta} \in R_{J,+}$, $\vec{\beta} \ne \vec{\beta}_i$ and 
$\langle \vec{\beta} , 
\vec{\beta}_i \rangle \ne 0$, then $s_i \vec{\beta} \in R_{J,+}$ with
\begin{equation}\label{11.3}
\langle s_i \vec{\beta} ,
\vec{\beta}_i \rangle =
- \langle \vec{\beta} ,
\vec{\beta}_i \rangle.
\end{equation}
Thus we can rewrite that product as
$$
O_J^{(2)} = \prod_{\vec{\beta} \in R_{J,+} \atop
\vec{\beta} = \vec{e}_{i+1} - \vec{e}_k \, (k > i+1)}
\Big ( D_{\vec{\beta}} + {1 \over \alpha} \Big )
\Big ( D_{s_i \vec{\beta}} + {1 \over \alpha} \Big ) \  .
\prod_{\vec{\beta} \in R_{J,+} \atop
\vec{\beta} = \vec{e}_{k} - \vec{e}_i \, (k < i)}
\Big ( D_{\vec{\beta}} + {1 \over \alpha} \Big )
\Big ( D_{s_i \vec{\beta}} + {1 \over \alpha} \Big ),
$$
and then use (\ref{11.1}) and (\ref{11.3}) to deduce
\begin{eqnarray}
&& s_i \prod_{\vec{\beta} \in R_{J,+} \atop
\vec{\beta} = \vec{e}_{i+1} - \vec{e}_k \, (k > i+1)}
\Big ( D_{\vec{\beta}} + {1 \over \alpha} \Big ) 
\Big ( D_{s_i \vec{\beta}}  + {1 \over \alpha} \Big )\nonumber \\
&=& \prod_{\vec{\beta} \in R_{J,+} \atop
\vec{e}_{i+1} - \vec{e}_k \, (k > i+1)}
\Big ( D_{s_i \vec{\beta}} + {1 \over \alpha} - {s_i \over \alpha}
\Big ) \Big ( D_{\vec{\beta}}  + {1 \over \alpha} +
{s_i \over \alpha} \Big ) s_i\nonumber \\
&=& \prod_{\vec{\beta} \in R_{J,+} \atop
\vec{e}_{i+1} - \vec{e}_k \, (k > i+1)}
\Big ( D_{s_i \vec{\beta}} + {1 \over \alpha} \Big ) 
\Big ( D_{\vec{\beta}}  + {1 \over \alpha} \Big ) s_i
\end{eqnarray}
where to obtain the second equality further use has been made of 
(\ref{11.1}).
Similarly, 
\begin{eqnarray}
&& s_i \prod_{\vec{\beta} \in R_{J,+} \atop
\vec{\beta} = \vec{e}_k - \vec{e}_i \, (k < i)}
\Big ( D_{\vec{\beta}} + {1 \over \alpha} \Big ) 
\Big ( D_{s_i \vec{\beta}}  + {1 \over \alpha} \Big )\nonumber \\
&=& \prod_{\vec{\beta} \in R_{J,+} \atop
\vec{e}_k - \vec{e}_i \, (k < i)}
\Big ( D_{s_i \vec{\beta}} + {1 \over \alpha} - {s_i \over \alpha}
\Big ) \Big ( D_{\vec{\beta}}  + {1 \over \alpha} +
{s_i \over \alpha} \Big ) s_i\nonumber \\
&=& \prod_{\vec{\beta} \in R_{J,+} \atop
\vec{e}_k - \vec{e}_i \, (k < i)}
\Big ( D_{s_i \vec{\beta}} + {1 \over \alpha} \Big ) 
\Big ( D_{\vec{\beta}}  + {1 \over \alpha} \Big ) s_i,
\end{eqnarray}
and so 
\begin{equation}\label{O2}
s_i O_J^{(2)} = O_J^{(2)} s_i.
\end{equation}

In the case that $\vec{\beta}_i \in R_{J,+}$, or equivalently
$i \in J$, we must also consider the action of $s_i$ on the final
factor in (\ref{11.2}). For this we see from (\ref{11.1}) and
(\ref{1.7a}) that
$$
s_i \Big ( D_{\vec{\beta}_i } + {1 \over \alpha} \Big ) =
- D_{\vec{\beta}_i } s_i + {2 \over \alpha} + {s_i \over \alpha},
$$
and thus for any $f$ such that $s_i f = -f$,
\begin{equation}\label{O3}
s_i \Big ( D_{\vec{\beta}_i } + {1 \over \alpha} \Big ) f =
\Big ( D_{\vec{\beta}_i } + {1 \over \alpha} \Big ) f.
\end{equation}

An immediate consequence of (\ref{O1}), (\ref{O2}) and
(\ref{O3}) is the action of $s_i$ on
$O_J S_{\eta^*}^{(I,J)}(z)$, and we can deduce from this that
the latter is proportional to $ S_{\eta^*}^{(I \cup J, \emptyset)}(z)$. 

\begin{prop}
For all $i \in I \cup J$,
\begin{equation}\label{O3'}
s_i \Big ( O_J S_{\eta^*}^{(I,J)}(z) \Big ) =
O_J S_{\eta^*}^{(I,J)}(z).
\end{equation}
Moreover, with
\begin{equation}\label{12.8}
c_{\eta^*} := \prod_{\vec{\beta} \in R_{J,+}}
\Big ((\overline{\eta^*}_{\vec{\beta}} + 1  )/ \alpha \Big )
\end{equation}
where $ \overline{\eta^*}_{\beta} :=  \overline{\eta^*}_j -  
\overline{\eta^*}_k$ for
$\vec{\beta}$ as in (\ref{1.5}), we have
\begin{equation}\label{12.9}
O_J S_{\eta^*}^{(I,J)}(z) = c_{\eta^*} S_{\eta^*}^{(I\cup J,\emptyset)}(z).
\end{equation}
\end{prop}

\noindent
Proof. \quad As already remarked, (\ref{O3'}) is an immediate
consequence of  (\ref{O1}), (\ref{O2}) and
(\ref{O3}). Also, from (\ref{11.0}), (\ref{5.2}) and (\ref{1.7b})
\begin{equation}\label{O3e}
O_J S_{\eta^*}^{(I,J)}(z) = \sum_{\mu \in W_{I \cup J}(\eta^*)}
\hat{c}_{\eta^* \mu} c_\mu E_\mu(z),
\end{equation}
where $c_\mu$ is specified by (\ref{12.8}). We know from
the remark at the end of the paragraph containing
(\ref{4.3}) that this structural formula together
with (\ref{O3'}) imply (\ref{12.9}).
\hfill $\square$

\medskip
Next we turn our attention to the structure of the image of the
product
$$
\Big ( \prod_{\vec{\beta} \in R_{J,+}} (z^{\vec{\beta}} - 1) \Big )
U_{\eta^*}^{(I,J)}(z)
$$
(c.f.~(\ref{2.I.1})) under the action of $O_J$.

\begin{prop}\label{p13}
Let $F(z)$ be an analytic function of $z_1,\dots,z_N$ in the
neighbourhood of $z=1^N$, and let $\Phi \subseteq R_+$.
Suppose $0 \le l \le \#(\Phi)$, then for any $\vec{\beta}_1,\dots, 
\vec{\beta}_l$,
$$
D_{\vec{\beta}_1} \cdots D_{\vec{\beta}_l}
\Big (  \prod_{\vec{\beta} \in \Phi}
(z^{\vec{\beta}} - 1) \, F(z) \Big ) =
\sum_{\Omega \subseteq R_+ \atop
\#(\Omega) = \#(\Phi) - l} 
h_\Omega(z) \prod_{\vec{\beta} \in \Omega}
(z^{\vec{\beta}} - 1)
$$
where $h_\Omega(z)$ is an analytic function in the neighbourhood of $z=1^N$ 
dependent on $\Omega$, $F$, and $\vec{\beta}_1,\dots,\vec{\beta}_l$.
\end{prop}

\noindent
Proof. \quad It suffices to establish the result for $l=1$, as the
general $l$ result follows by induction. Now, for general $f$ and $g$
and $\vec{\beta} = \vec{e}_j - \vec{e}_k$ $(j < k)$ we have
\begin{equation}\label{gz}
{1 \over 1 - z^{\vec{\beta}}}(1 - s_{\vec{\beta}}) (fg) =
{f - s_{\vec{\beta}} f \over 1 - z^{\vec{\beta}} } g
+ {s_{\vec{\beta}} f \over 1 - z^{\vec{\beta}} }
(1 -  s_{\vec{\beta}} ) g. 
\end{equation}
We apply this formula with
\begin{equation}\label{gza}
f = F \prod_{\vec{\gamma} \in \Phi \atop (\vec{\gamma}, \vec{\beta})
= 0} ( z^{\vec{\gamma}} - 1), \qquad
g = \prod_{\vec{\gamma} \in \Phi \atop (\vec{\gamma}, \vec{\beta}) \ne 0}
( z^{\vec{\gamma}} - 1).
\end{equation}
Since 
$$
f - s_{\vec{\beta}} f =
\prod_{\vec{\gamma} \in \Phi \atop (\vec{\gamma}, \vec{\beta}) = 0} 
(z^{\vec{\gamma}} - 1) (F - s_{\vec{\beta}}F),
$$
it follows immediately that the first term has the required structure.

To show that the second term has the required structure, define a total
order $<_{\vec{\beta}}$ on $\{\vec{\gamma} \in \Phi |
(\vec{\gamma}, \vec{\beta}) \ne 0 \}$ by the requirement that
$\vec{\gamma}$ and $s_{\vec{\beta}}(\vec{\gamma})$ are adjacent in the order.
Otherwise the order is arbitrary. In terms of this order we can write
\begin{equation}\label{gz1}
(1 - s_{\vec{\beta}}) \prod_{\vec{\gamma} \in \Phi \atop 
(\vec{\gamma}, \vec{\beta}) \ne 0} ( z^{\vec{\gamma}} - 1)
=
\sum_{\vec{\alpha} \in \Phi}
\prod_{\vec{\gamma} \in \Phi \atop
\vec{\gamma} <_{\vec{\beta}} \vec{\alpha} }
(z^{s_{\vec{\beta}}(\vec{\gamma})} - 1) \Big (
(z^{\vec{\alpha}} - 1) -
(z^{s_{\vec{\beta}}(\vec{\alpha})} - 1) \Big )
\prod_{\vec{\gamma} \in \Phi \atop
\vec{\alpha} <_{\vec{\beta}} \vec{\gamma} }
(z^{\vec{\gamma}} - 1),
\end{equation}
by the telescoping of the sum (the same mechanism responsible for
(\ref{gz})). Now for a given $\vec{\alpha}$, if vectors $\vec{\gamma}_1$ 
and $\vec{\gamma}_2$ satisfy $\vec{\gamma}_1 <_{\vec{\beta}} \vec{\alpha} 
<_{\vec{\beta}} \vec{\gamma}_2$, then $\vec{\gamma}_1$ and $\vec{\gamma}_2$ 
are not adjacent so that $s_{\vec{\beta}}(\vec{\gamma}_1) \ne \vec{\gamma}_2$.
It follows that all factors in (\ref{gz1}) are distinct. Because
$(z^{\vec{\alpha}} - z^{s_{\vec{\beta}}(\vec{\alpha})})$ is divisible
by $(1 - z^{\vec{\beta}})$, we see that the second term in
(\ref{gz}) with the substitution (\ref{gza}) has the sought structure.
\hfill $\square$

\medskip
Proposition \ref{p13} can be used to establish the following
structural formula.

\begin{prop}\label{p14}
Let $G(z)$ be an analytic function of $z_1,\dots,z_N$ in the
neighbourhood of $z=1^N$, and let $F(z)$, $\Phi$ and the
$\vec{\beta}_i$ be as in Proposition \ref{p13}.  Then for 
$0 \le l \le \#(\Phi)$, there exist functions $\tilde{h}_\Omega(z)$
analytic in a neighbourhood of $z = 1^N$ such that
\begin{eqnarray}\label{15.1}
&&D_{\vec{\beta}_1} \cdots D_{\vec{\beta}_l}
\Big (  \prod_{\vec{\beta} \in \Phi}
(z^{\vec{\beta}} - 1) \, F(z) G(z)\Big ) -
G(z) D_{\vec{\beta}_1} \cdots D_{\vec{\beta}_l}
\Big (  \prod_{\vec{\beta} \in \Phi}
(z^{\vec{\beta}} - 1) \, F(z) \Big ) \nonumber \\ &&
=
\sum_{\Omega \subseteq R_+ \atop
\#(\Omega) = \#(\Phi) - l+1} 
\tilde{h}_\Omega(z) \prod_{\vec{\beta} \in \Omega}
(z^{\vec{\beta}} - 1)
\end{eqnarray}
\end{prop}

\noindent
Proof. \quad 
Consider first the case $l=1$. We can check from the
definition (\ref{1.2}) that
\begin{eqnarray}\label{15.2}
&&D_{\vec{\lambda}} \Big (  \prod_{\vec{\beta} \in \Phi}
(z^{\vec{\beta}} - 1) \, F(z) G(z)\Big ) \nonumber \\&& \quad =
G(z) D_{\vec{\lambda}}
\Big (  \prod_{\vec{\beta} \in \Phi}
(z^{\vec{\beta}} - 1) \, F(z) \Big )  
+ \prod_{\vec{\beta} \in \Phi}
(z^{\vec{\beta}} - 1) \, F(z)
\Big ( \sum_{j=1}^N \lambda_j z_j {\partial \over \partial z_j}
\Big ) G(z)  \nonumber \\&& \qquad  + {1 \over \alpha}
\sum_{\vec{\beta} \in R_+}
{\langle \vec{\lambda}, \vec{\beta} \rangle \over
z^{\vec{\beta}} - 1 }
\Big ( G(z) - s_{\vec{\beta}} G(z) \Big )
s_{\vec{\beta}} \Big (
\prod_{\vec{\gamma} \in \Phi}
(z^{\vec{\gamma}} - 1) \, F(z) \Big ).
\end{eqnarray}
The structure (\ref{15.1}) follows immediately. We now proceed
inductively, assuming the result (\ref{15.1}) for some
$1 \le l < \#(\Phi)$, with our task being to then prove its
validity for $l \mapsto l+1$.  Since $\vec{\beta}_1,\dots,\vec{\beta}_l$
are arbitrary in the statement of the proposition, then for $\vec{\beta}_1,
\vec{\beta}_2, \dots,\vec{\beta}_{l+1}$,
\begin{eqnarray}\label{15.3}
&&D_{\vec{\beta}_2} \cdots D_{\vec{\beta}_{l+1}}
\Big (  \prod_{\vec{\beta} \in \Phi}
(z^{\vec{\beta}} - 1) \, F(z) G(z)\Big ) -
G(z) D_{\vec{\beta}_2} \cdots D_{\vec{\beta}_{l+1}}
\Big (  \prod_{\vec{\beta} \in \Phi}
(z^{\vec{\beta}} - 1) \, F(z) \Big ) \nonumber \\ &&
=
\sum_{\Omega \subseteq R_+ \atop
\#(\Omega) = \#(\Phi) - l+1} 
\tilde{h}_\Omega(z) \prod_{\vec{\beta} \in \Omega}
(z^{\vec{\beta}} - 1)
\end{eqnarray}
for some $\tilde{h}_\Omega(z)$ which are analytic in a neighbourhood 
of $z = 1^N$.
Applying $D_{\vec{\beta}_1}$ to
both sides of (\ref{15.3}) and making use of Proposition \ref{p13}
on the RHS we see that
\begin{eqnarray}\label{16.1}
&&D_{\vec{\beta}_1} \cdots D_{\vec{\beta}_{l+1}}
\Big (  \prod_{\vec{\beta} \in R_\Phi}
(z^{\vec{\beta}} - 1) \, F(z) G(z)\Big ) -
D_{\vec{\beta}_1}
\Big \{ G(z) D_{\vec{\beta}_2} \cdots D_{\vec{\beta}_{l+1}}
\Big (  \prod_{\vec{\beta} \in R_\Phi}
(z^{\vec{\beta}} - 1) \, F(z) \Big ) \Big \}  \nonumber \\&& \qquad
=
\sum_{\Omega \subseteq R_+ \atop
\#(\Omega) = \#(\Phi) - l}
{h}_\Omega(z) \prod_{\vec{\beta} \in \Omega}
(z^{\vec{\beta}} - 1)
\end{eqnarray}
for some $h_\Omega(z)$ which are analytic in a neighbourhood of $z=1^N$.
Now Proposition \ref{p13} also gives that
$$
D_{\vec{\beta}_2} \cdots D_{\vec{\beta}_{l+1}}
\Big (  \prod_{\vec{\beta} \in \Phi}
(z^{\vec{\beta}} - 1) \, F(z) \Big ) =
\sum_{\Omega \subseteq R_+ \atop
\#(\Omega) = \#(\Phi) - l}
\hat{h}_\Omega(z) \prod_{\vec{\beta} \in \Omega}
(z^{\vec{\beta}} - 1)
$$
for some $\hat{h}_\Omega(z)$ which are analytic in a neighbourhood of $z=1^N$.
Application of the already established
$l=1$ case of the present proposition then shows
\begin{eqnarray}
&&D_{\vec{\beta}_1}
\Big \{ G(z) D_{\vec{\beta}_2} \cdots D_{\vec{\beta}_{l+1}}
\Big (  \prod_{\vec{\beta} \in \Phi}
(z^{\vec{\beta}} - 1) \, F(z) \Big ) \Big \} -
G(z) D_{\vec{\beta}_1} \cdots D_{\vec{\beta}_{l+1}}
\Big (  \prod_{\vec{\beta} \in R_\Phi}
(z^{\vec{\beta}} - 1) \, F(z) \Big )  \nonumber \\
&& =
\sum_{\Omega \subseteq R_+ \atop
\#(\Omega) = \#(\Phi) - l}
{h}'_\Omega(z) \prod_{\vec{\beta} \in \Omega}
(z^{\vec{\beta}} - 1),
\end{eqnarray}
for some ${h}'_\Omega(z)$ which are analytic in a neighbourhood of $z=1^N$,
and this substituted in (\ref{16.1})
establishes (\ref{15.1}) in the case $l \mapsto l + 1$.
\hfill $\square$

\medskip
An immediate corollary of Proposition \ref{p14} is the following
evaluation identity.

\begin{cor}
In the notation of Proposition \ref{p13}, for $0 \le l \le \#\Phi - 1$,
\begin{equation}\label{17.1}
D_{\vec{\beta}_1} \cdots D_{\vec{\beta}_l}
\Big (  \prod_{\vec{\beta} \in \Phi}
(z^{\vec{\beta}} - 1) \, F(z) \Big ) \Big |_{z = 1^N} = 0
\end{equation}
while
\begin{equation}\label{17.2}
D_{\vec{\beta}_1} \cdots D_{\vec{\beta}_{\# \Phi}} 
\Big (  \prod_{\vec{\beta} \in \Phi}
(z^{\vec{\beta}} - 1) \, F(z) \Big ) \Big |_{z = 1^N} =
F(z) \Big |_{z = 1^N}
D_{\vec{\beta}_1} \cdots D_{\vec{\beta}_{\# \Phi}}
\Big (  \prod_{\vec{\beta} \in \Phi}
(z^{\vec{\beta}} - 1)  \Big ) \Big |_{z = 1^N}.
\end{equation}
In particular
\begin{equation}\label{17.2a}
O_J S_{\eta^*}^{(I,J)}(z) \Big |_{z = 1^N} =
U_{\eta^*}^{(I,J)}(1^N)
\Big ( \prod_{\vec{\beta} \in R_{J,+}} D_{\vec{\beta}} \Big )
\Big (  \prod_{\vec{\beta} \in R_{J,+}} (
z^{\vec{\beta}} - 1)  \Big ) \Big |_{z = 1^N}.
\end{equation}
\end{cor}

\noindent
Proof. \quad The equations (\ref{17.1}) and (\ref{17.2}) are
immediate consequences of Proposition \ref{p14}. Using these
equations, (\ref{17.2a}) follows after recalling the definitions
(\ref{2.I.1}) and (\ref{11.0}), and noting that
\begin{equation}
{\prod_{\vec{\beta} \in R_{J,+}} (
z^{\vec{\beta}} - 1) \over \prod_s \Delta_{\tilde{J}_s}(z) } 
= z^{\vec{\lambda}}
\end{equation}
for some $\vec{\lambda}$.
\hfill $\square$

Substituting (\ref{12.9}), (\ref{12.8}) and Proposition \ref{p10} in
(\ref{17.2a}) we see that
\begin{equation}\label{17.3}
U_{\eta^*}^{(I,J)}(1^N) =
{|{\rm Sym} \, I \cup J | \over {\cal M}_{ I , \eta^*}}
{e_{\eta^+} \over d_{\eta^{*(-,-)}}} {1 \over k_{J}}
 \prod_{\vec{\beta} \in R_{J,+}} \Big ( (\overline{\eta^*}_{\vec{\beta}} + 1)/
\alpha \Big ) 
\end{equation}
where
\begin{equation}\label{kJ}
k_{J} := \Big ( \prod_{\vec{\beta} \in R_{J,+}} D_{\vec{\beta}} \Big )
\Big (  \prod_{\vec{\beta} \in R_{J,+}} (
z^{\vec{\beta}} - 1)  \Big ) \Big |_{z = 1^N}.
\end{equation}
With our objective being to compute $U_{\eta^*}^{(I,J)}(1^N)$, we see
from (\ref{17.3}) that the remaining task is to compute $k_{J}$
as specified by (\ref{kJ}).

Suppose for a given $J$ we could choose a composition $\eta^*$ such
that $U_{\eta^*}(1^N)$ could be evaluated directly from its definition
(\ref{2.I.1}). Then because $k_{J}$ is independent of $\eta^*$, that
evaluation
substituted in (\ref{17.3}) will specify $k_{J}$. To implement 
this strategy, we consider the particular composition
$\eta^* = \delta$  defined by
\begin{equation}\label{18.1}
\delta_j = \left \{ \begin{array}{ll}
\max(\tilde{J}_s)  - j, & j \in \tilde{J}_s \\
0, & j \notin \tilde{J} \end{array} \right.
\end{equation} 
With this definition and the final equality in (\ref{3.2}) we see that
$$
{\rm Asym}_J \, z^{\delta} = \Delta_{\tilde{J}}(z)
$$
and thus according to (\ref{4.3})
$$
S_{\delta}^{(\emptyset, J)}(z) = \Delta_{\tilde{J}}(z).
$$
Substituting this in (\ref{2.I.1}) shows
\begin{equation}\label{18.1a}
U_{\delta}^{(\emptyset, J)}(1^N) = 1.
\end{equation}
The value of the quantity $k_J$ can now be deduced by substituting
(\ref{18.1a}) in (\ref{17.3}) and simplifying the result. 

\begin{prop}
With $k_{J}$ specified by (\ref{kJ}) we have
\begin{equation}\label{kJ1}
k_{J} = \prod_{s} |\tilde{J}_s|! \, \alpha^{-|\delta|} e_{\delta}.
\end{equation}
\end{prop}

\noindent
Proof. \quad
Substituting (\ref{18.1}) in (\ref{17.3}) gives
$$
k_{J} = \prod_{s} |\tilde{J}_s|! \, \alpha^{-|\delta|} {e_{\delta^+} \over
d_{\delta^{(\cdot,-)}}}
\prod_{\vec{\beta} \in R_{J,+}}(\bar{\delta}_{\vec{\beta}} + 1) 
$$
so to obtain (\ref{18.1}) we must show
\begin{equation}\label{kJ2}
d_{\delta^{(\cdot,-)}} = \prod_{\vec{\beta} \in R_{J,+}}
(\bar{\delta}_{\vec{\beta}} + 1) =
\prod_{(i,j) \in \eta^*}(
\bar{\delta}_i -
\bar{\delta}_{i+j} + 1),
\end{equation}
where the final equality follows from the definitions of the
quantities in the preceeding expression. Now after recalling
the definition (\ref{5.3}) of $d_\eta$ and the meaning of
$\delta^{(\cdot,-)}$, we see the LHS of
(\ref{kJ2}) can be rewritten
\begin{equation}\label{kJ3}
\prod_s \prod_{i'=0}^{J_s^{\rm max} - J_s^{\rm min}}
\prod_{j=1}^{i'+1} \Big (
\alpha( a(J_s^{\rm min} + i' + 1, j) + 1) + l(J_s^{\rm min} + i' + 1,j) + 1 
\Big )
\end{equation}
where the arm and leg lengths are to be calculated with respect to the
diagram of $\delta^{(\cdot,-)}$. Similarly, on the RHS of
(\ref{kJ2}) we can write
\begin{equation}\label{kJ4}
\prod_s \prod_{i'=0}^{J_s^{\rm max} - J_s^{\rm min}}
\prod_{j=1}^{i'+1} \Big (
\bar{\delta}_{J_s^{\rm max} - i'}  -
\bar{\delta}_{J_s^{\rm max} - i' + j} +1 \Big ).
\end{equation}
It is easy to check from the definitions that the terms in the
products (\ref{kJ3}) and (\ref{kJ4}) are in fact identical, and
thus (\ref{kJ2}) does indeed hold.
\hfill $\square$

Substituting  (\ref{kJ1}) in (\ref{17.3}) 
and noting $|{\rm Sym} \, I \cup J| =
(\prod_s |\tilde{I}_s|!) (\prod_s |\tilde{J}_s|!)$
gives the sought evaluation
formula.
\begin{prop}\label{p20}
We have
\begin{equation}\label{20.3}
U_{\eta^*}^{(I,J)}(1^N) =
{\prod_s | \tilde{I}_s|!
 \over   {\cal M}_{ I , \eta^*}}
{e_{\eta^+}  \over \alpha^{-|\delta|}
e_{\delta^+}d_{\eta^{*(-,-)}}} 
 \prod_{\vec{\beta} \in R_{J,+}} \Big ( (\overline{\eta^*}_{\vec{\beta}} + 1)/
\alpha\Big )
\end{equation}
\end{prop}

\section{Discussion}
\setcounter{equation}{0}
Let us first compare our result (\ref{20.3}) with that which can be deduced
from the evaluation formula of Dunkl \cite{Du98}.
In our notation, the latter formula is
\begin{equation}\label{dn1}
{{\cal O}_{I,J}E_\eta(z) \over
\prod_s \Delta_{\tilde{J}_s}(z)} \bigg |_{z = 1^N} =
\prod_s | \tilde{I}_s|!
{ e_{\eta^+} \over \alpha^{-|\delta|} e_{\delta^+}}
{1 \over d_{\eta^*}}
\prod_{\vec{\beta} \in R_{J,+}} \Big (( \overline{\eta^*}_{\vec{\beta}} - 1 )/
\alpha \Big ).
\end{equation}
Making use of (\ref{4.3}) and (\ref{8.1}) in the case $\eta = \eta^*$,
and recalling (\ref{2.I.1}),
we see that this implies
\begin{equation}\label{dn2}
U_{\eta^*}^{(I,J)}(1^N) =
{\prod_s | \tilde{I}_s|! \over {\cal M}_{I,\eta^*}}
{d_{\eta^{*(-,-)}}' \over
d_{\eta^{*(-,0)}}' d_{\eta^{*(-,0)}} }
{e_{\eta^+} \over \alpha^{-|\delta|} e_{\delta^+}}
\prod_{\vec{\beta} \in R_{J,+}} \Big (( \overline{\eta^*}_{\vec{\beta}} - 1 )/
\alpha \Big ).
\end{equation}
Note that this
is a different form to the result (\ref{20.3}). To reconcile the
two forms we must have 
\begin{equation}\label{dn4}
{d_{\eta^{*(-,-)}}' d_{\eta^{*(-,-)}} \over
d_{\eta^{*(-,0)}}' d_{\eta^{*(-,0)}}}
= \prod_{\vec{\beta} \in R_{J,+}} {\overline{\eta^*}_{\vec{\beta}} + 1
\over \overline{\eta^*}_{\vec{\beta}} - 1}.
\end{equation}
In fact (\ref{dn4}) is just a special case of the following result.

\begin{prop}
For $w \in W_J$ we have
\begin{equation}\label{dn5}
{ d_{w \eta^{*(-,0)}}' d_{w \eta^{*(-,0)}} \over
d_{\eta^{*(-,0)}}' d_{\eta^{*(-,0)}} }
= \prod_{\vec{\beta} \in R_{J,+} \cap w^{-1} R_{J,-}} 
{\overline{\eta^*}_{\vec{\beta}} + 1
\over \overline{\eta^*}_{\vec{\beta}} - 1}.
\end{equation}
\end{prop}

\noindent Proof. \quad
This is demonstrated by induction of $l(w)$ 
(the length of $w$, i.e. the smallest value of $l$ for which $w$ can 
be expressed as a product $w = s_{j_1} \dots s_{j_l}$ for $j_1,\dots,j_l 
\in J$).
First, note that for all $i$, and for all compositions $\mu$ such that 
$s_i \mu \ne \mu$, $\overline{s_i \mu} = s_i \overline{\mu}$.
Since $I \cap J = \emptyset$, and $I$ and $J$ satisfy (\ref{es1}), then 
it follows that for all $i \in I$ and all compositions $\mu$, 
$\overline{s_i \mu}_j = \bar{\mu}_j$ for all $j \in \tilde{J}$.  
This generalizes to the result that if $w \in W_I$, then for all compositions 
$\mu$, $\overline{w \mu}_j = \bar{\mu}_j$ for $j \in \tilde{J}$.
A consequence of this result is that $\overline{\eta^{*(-,0)}}_j 
= \overline{\eta^*}_j$ for $j \in \tilde{J}$, and 
$\overline{\eta^{*(-,0)}}_{\vec{\beta}} = \overline{\eta^*}_{\vec{\beta}}$\ \ 
for $\vec{\beta} \in R_{J,+}$.
Now, (\ref{dn5}) is clearly true when $w = 1$ (the left hand side is 1, 
and the product on the right hand side is empty, and therefore evaluates 
to 1).
For general $w \in W_J$, $w = s_{j_1} \dots s_{j_l}$ for some $j_1,\dots,
j_l \in J$, where $l = l(w)$.
It follows that $w^{-1} \vec{\alpha}_{j_1} \in R_{J,-}$, where 
$\vec{\alpha}_j = \vec{e}_j - \vec{e}_{j+1}$ for $j = 1,\dots,N-1$.
Let $w_1 = s_{j_2} \dots s_{j_l}$, then $l(w_1) = l-1$, $w = s_{j_1} w_1$,
and $w_1^{-1} \vec{\alpha}_{j_1} \in R_{J,+}$.
Since $w_1^{-1} \vec{\alpha}_{j_1} \in R_{J,+}$, then 
$(w_1 \eta^{*(-,0)})_{j_1} > (w_1 \eta^{*(-,0)})_{j_1+1}$, and so
$$
{ d_{w \eta^{*(-,0)}}' \over d_{w_1 \eta^{*(-,0)}}' } 
= { \overline{w_1 \eta^{*(-,0)}}_{\vec{\alpha}_{j_1}} \over 
\overline{w_1 \eta^{*(-,0)}}_{\vec{\alpha}_{j_1}} - 1 },
$$
and 
$$
{ d_{w \eta^{*(-,0)}} \over d_{w_1 \eta^{*(-,0)}} } 
= { \overline{w_1 \eta^{*(-,0)}}_{\vec{\alpha}_{j_1}} + 1 \over 
\overline{w_1 \eta^{*(-,0)}}_{\vec{\alpha}_{j_1}} }.
$$
It follows that 
$$
{ d_{w \eta^{*(-,0)}}' d_{w \eta^{*(-,0)}} \over 
d_{w_1 \eta^{*(-,0)}}' d_{w_1 \eta^{*(-,0)}} } 
= { \overline{w_1 \eta^{*(-,0)}}_{\vec{\alpha}_{j_1}} + 1 \over 
\overline{w_1 \eta^{*(-,0)}}_{\vec{\alpha}_{j_1}} - 1 }.
$$
As a consequence of (\ref{es}), then $\overline{w_1 \eta^{*(-,0)}} 
= w_1 \overline{\eta^{*(-,0)}}$, so that 
\begin{eqnarray}
{ d_{w \eta^{*(-,0)}}' d_{w \eta^{*(-,0)}} \over 
d_{w_1 \eta^{*(-,0)}}' d_{w_1 \eta^{*(-,0)}} } 
&=& { \overline{\eta^{*(-,0)}}_{w_1^{-1} \vec{\alpha}_{j_1}} + 1 \over 
\overline{\eta^{*(-,0)}}_{w_1^{-1} \vec{\alpha}_{j_1}} - 1 }\\
&=& { \overline{\eta^*}_{w_1^{-1} \vec{\alpha}_{j_1}} + 1 \over 
\overline{\eta^*}_{w_1^{-1} \vec{\alpha}_{j_1}} - 1 }.
\end{eqnarray}
By the induction hypothesis, 
$$
{ d_{w_1 \eta^{*(-,0)}}' d_{w_1 \eta^{*(-,0)}} \over
d_{\eta^{*(-,0)}}' d_{\eta^{*(-,0)}} }
= \prod_{\vec{\beta} \in R_{J,+} \cap w_1^{-1} R_{J,-}} 
{ \overline{\eta^*}_{\vec{\beta}} + 1
\over \overline{\eta^*}_{\vec{\beta}} - 1 },
$$
since $l(w_1) = l-1$, and so 
\begin{equation}\label{dn6}
{ d_{w \eta^{*(-,0)}}' d_{w \eta^{*(-,0)}} \over
d_{\eta^{*(-,0)}}' d_{\eta^{*(-,0)}} }
= { \overline{\eta^*}_{w_1^{-1} \vec{\alpha}_{j_1}} + 1 \over 
\overline{\eta^*}_{w_1^{-1} \vec{\alpha}_{j_1}} - 1 } 
\prod_{\vec{\beta} \in R_{J,+} \cap w_1^{-1} R_{J,-}} 
{ \overline{\eta^*}_{\vec{\beta}} + 1 \over 
\overline{\eta^*}_{\vec{\beta}} - 1 }.
\end{equation}
Suppose $\vec{\beta} \in R_{J,+} \cap w^{-1} R_{J,-}$, then $w \vec{\beta} 
= s_{j_1} w_1 \vec{\beta} \in R_{J,-}$, and so $w_1 \vec{\beta} \in R_{J,-}$ 
or $w_1 \vec{\beta} = \vec{\alpha}_{j_1}$.  In the first case, 
$\vec{\beta} \in R_{J,+} \cap w_1^{-1} R_{J,-}$, and in the second, 
$\vec{\beta} = w_1^{-1} \vec{\alpha}_{j_1}$.
Conversely, if $\vec{\beta} \in R_{J,+} \cap w_1^{-1} R_{J,-}$, then 
$w_1 \vec{\beta} \in R_{J,-}$, and $w_1 \vec{\beta} \ne - \vec{\alpha}_{j_1}$
since $w_1^{-1} \vec{\alpha}_{j_1} \in R_{J,+}$, and so 
$w \vec{\beta} = s_{j_1} w_1 \vec{\beta} \in R_{J,-}$.
Alternatively, if $\vec{\beta} = w_1^{-1} \vec{\alpha}_{j_1}$, then 
$\vec{\beta} \in R_{J,+}$ and $w \vec{\beta} = s_{j_1} w_1 \vec{\beta} 
= - \vec{\alpha}_{j_1} \in R_{J,-}$.
Since $\vec{\beta} \in R_{J,+} \cap w^{-1} R_{J,-}$ is equivalent to 
$\vec{\beta} \in R_{J,+} \cap w_1^{-1} R_{J,-}$ or $\vec{\beta} 
= w_1^{-1} \vec{\alpha}_{j_1}$, then (\ref{dn6}) becomes
$$
{ d_{w \eta^{*(-,0)}}' d_{w \eta^{*(-,0)}} \over
d_{\eta^{*(-,0)}}' d_{\eta^{*(-,0)}} }
= \prod_{\vec{\beta} \in R_{J,+} \cap w^{-1} R_{J,-}} 
{ \overline{\eta^*}_{\vec{\beta}} + 1 \over 
\overline{\eta^*}_{\vec{\beta}} - 1 },
$$
which is just (\ref{dn5}) for $w$.
This completes the proof by induction.
\hfill $\square$

To deduce (\ref{dn4}) from (\ref{dn5}),
let $w_{0,J}$ be the longest element of $W_J$, so that $w_{0,J} \eta^{*(-,0)} 
= \eta^{*(-,-)}$ and $w_{0,J} R_{J,+} = R_{J,-}$. Then, substituting 
$w = w_{0,J}$ into (\ref{dn5}) gives
$$
{d_{\eta^{*(-,-)}}' d_{\eta^{*(-,-)}} \over
d_{\eta^{*(-,0)}}' d_{\eta^{*(-,0)}}}
= \prod_{\vec{\beta} \in R_{J,+} \cap w_{0,J}^{-1} R_{J,-}} 
{\overline{\eta^*}_{\vec{\beta}} + 1 \over 
\overline{\eta^*}_{\vec{\beta}} - 1}
= \prod_{\vec{\beta} \in R_{J,+}}
{\overline{\eta^*}_{\vec{\beta}} + 1 \over 
\overline{\eta^*}_{\vec{\beta}} - 1}
$$
as required.

Next we will contrast our strategy to derive (\ref{20.3}) with that
used by Dunkl to derive (\ref{dn2}). In our method the key ingredient
is the operator $O_J$ with its action (\ref{12.9}) and evaluation
property (\ref{17.2a}). In Dunkl's approach the key ingredient is the
skew operator $\psi_J$, constructed so that
\begin{equation}\label{du5}
\psi_J s_j = - s_j \psi_J, \qquad j \in J
\end{equation}
and possessing an evaluation property analogous to that of
$O_J$.
It follows from (\ref{du5}) that for any $f$ such that $s_j f = - f$,
\begin{equation}\label{du6}
s_j ( \psi_J f) = \psi_J f.
\end{equation}
The operator $O_J$ has the same action when restricted to this class of
$f$. However $O_J$ does not exhibit the general algebraic property
(\ref{du5}). Thus the two approaches are identical in strategy except
that in Dunkl's work a more complicated operator $\psi_J$ is
constructed which plays the role of our $O_J$.

The result (\ref{dn1}), or equivalently (\ref{20.3}), has been used by
Kato and Yamamoto \cite{KY98} to obtain the expansion of 
(\ref{I.3.1}) in terms of $\{U_{\eta^*}^{(I,J)}\}$
(Proposition \ref{pKY} below). This result was then used in the exact
computation of the retarded Green function for spin generalizations of
the Calogero-Sutherland system (\ref{I.1}). Here we will use the
expansion formula to derive a constant term identity, which
generalizes an integration formula due to Macdonald, Kadell and
Kaneko \cite{Ma87, Ka97g, Ka93}. Let us first present the expansion
formula.

\begin{prop}\label{pKY}
For general $r$ we have
\begin{equation}\label{N.0}
\prod_{i \notin \tilde{J}}(1 - x_i)^{r-1}
\prod_s \prod_{j \in \tilde{J}_s} (1 - x_j)^{r - | \tilde{J}_s|} =
\sum_{\eta^*} { \alpha^{|\eta|} \over d_{\eta^*}'}
{ [1 - r]_{\eta^+}^{(\alpha)} \over
[1 - r]_{\delta^+}^{(\alpha)} }
\prod_{\vec{\beta} \in R_{J,+}} \Big ((\overline{\eta^*}_{\vec{\beta}} - 1)/
\alpha \Big ) U_{\eta^*}^{(I,J)}(x).
\end{equation}
\end{prop}

\noindent
Proof. \quad We recall the nonsymmetric Cauchy product expansion
\begin{equation}\label{N.1}
\Omega(x,y) := \prod_{i=1}^N {1 \over (1 - x_i y_i)}
\prod_{i,j=1}^N {1 \over (1 - x_i y_j)^{1/\alpha}} =
\sum_{\eta} {d_\eta \over d_\eta'} E_\eta(x) E_\eta(y)
\end{equation}
(c.f.~(\ref{t1})). Applying the operator ${\cal O}_{I,J}$ to the variables
$x$ on the RHS of (\ref{N.1}) gives
\begin{eqnarray}\label{N.2}
{\cal O}_{I,J} \Omega(x,y) & = & \sum_{\eta^*}
\sum_{\mu \in W_{I \cup J}(\eta^*)} {d_\mu \over d_\mu'}
a_\mu^{(I,J)} S_{\eta^*}^{(I,J)}(x) E_\mu(y) \nonumber \\
& = & \sum_{\eta^*} a_{\eta^*}^{-1} S_{\eta^*}^{(I,J)}(x)
S_{\eta^*}^{(I,J)}(y),
\end{eqnarray}
where the second equality follows from (\ref{9.1}), and
$a_{\eta^*}$ is given by (\ref{9.3}). On the other hand, applying the
operator ${\cal O}_{I,J}$ directly to the $x$ variables in its
definition as a product, and using the Cauchy formula
$$
{\rm Asym}^{(x)} \, \prod_{j \in \tilde{J}_s}
{1 \over (1 - x_j y_j)} =
{ \Delta_{\tilde{J}_s}(x) \Delta_{\tilde{J}_s}(y)  \over
\prod_{i,j \in \tilde{J}_s}(1 - x_i y_j)}
$$
we see that
\begin{eqnarray}\label{N.3}
&&
{\cal O}_{I,J} \Omega(x,y) \nonumber \\   = 
\Big ( \prod_{i \notin \tilde{I} \cup \tilde{J}}
{1 \over 1 - x_i y_i} \Big ) \Big ( {\rm Sym} \,
\prod_{i \in \tilde{I}} {1 \over 1 - x_i y_i} \Big )
\Big ( \prod_s { \Delta_{\tilde{J}_s}(x) \Delta_{\tilde{J}_s}(y)  \over
\prod_{i,j \in \tilde{J}_s}(1 - x_i y_j)} \Big ) 
\prod_{i,j=1}^N {1 \over (1 - x_i y_j)^{1/\alpha}}.
\end{eqnarray}
Equating (\ref{N.2}) and (\ref{N.3}), dividing by
$\prod_s \Delta_{\tilde{J}_s}(x) \Delta_{\tilde{J}_s}(y)$ and
substituting $y = 1^N$ we obtain
\begin{equation}\label{N.4}
\prod_s |\tilde{I}_s|! \prod_{i \notin \tilde{J}}
(1 - x_i)^{-1-N/\alpha} \prod_s \prod_{j \in \tilde{J}_s}
(1 - x_j)^{-|\tilde{J}_s| - N/\alpha} =
\sum_{\eta^*} a_{\eta^*}^{-1} U_{\eta^*}^{(I,J)}(x)
U_{\eta^*}^{(I,J)}(1^N).
\end{equation}

Now it follows from (\ref{9.3}) and (\ref{dn1}) that 
\begin{equation}\label{N.5}
{U_{\eta^*}^{(I,J)}(1^N) \over a_{\eta^*} \prod_s |\tilde{I}_s |!}
=  
{e_{\eta^+} \over \alpha^{-|\delta|} e_{\delta^+}}
{1 \over d_{\eta^*}'} \prod_{\vec{\beta} \in R_{J,+}}
\Big ( ( \overline{\eta^*}_{\vec{\beta}} - 1)/\alpha \Big ).
\end{equation}
The only factor in this expression dependent on the number of variables
$N$ is
\begin{equation}\label{N.6}
{e_{\eta^+} \over \alpha^{-|\delta|} e_{\delta^+}} =
\alpha^{|\eta|} {[1 + N/\alpha]_{\eta^+}^{(\alpha)} \over
[1 + N/\alpha]_{\delta^+}^{(\alpha)} },
\end{equation}
where to obtain the equality use has been made of (\ref{se}). At
order $|\eta|$ in $x$, both sides of (\ref{N.4}) are thus polynomials
in $N/\alpha$, which are equal for each $N=1,2,\dots$. Thus
$N/\alpha$ can be replaced by the continuous variable $-r$ and
(\ref{N.0}) follows. \hfill $\square$

Our constant term (CT) identity will be deduced from (\ref{N.0}) by
making use of the fact that $\{S_{\eta^*}^{(I,J)}(x)\}$ is
orthogonal with respect to the inner product
\begin{equation}\label{4.13}
\langle f, g \rangle_C := {\rm CT} \Big (
f(x) g(1/x) \Big ( \Delta(x) \Delta(1/x) \Big )^{1/\alpha}
\Big ),
\end{equation}
where $\Delta(x) := \prod_{1 \le j < k \le N} (x_k - x_j)$. This
constant term is well defined for $f$ and $g$ Laurent polynomials
and $1/\alpha \in \mathbb Z^+$. It is defined beyond these cases as a
multidimensional Fourier integral. We require the evaluation of the
norm of $S_{\eta^*}^{(I,J)}$ in this inner product \cite{BDF00,KY98}.

\begin{prop}
We have
\begin{equation}\label{M.8}
{|| S_{\eta^*}^{(I,J)} ||^2_C \over || 1 ||_C^2} =
{|{\cal O}_{I,J}| \over {\cal M}_{I,\eta^*}}
{d_{\eta^*}' d_{\eta^{*(-,-)}}' \over
d_{\eta^{*(-,0)}}' d_{\eta^{*(-,0)}}} {e_{\eta^+} \over
e_{\eta^+}'}
\end{equation}
where 
$e_{\eta^+}$ is as in (\ref{se}) while
\begin{equation}\label{4.14'}
e_{\eta^+}' = [1+(N-1)/\alpha]_{\eta^+}^{(\alpha)}
\end{equation}
\end{prop}

\noindent
Proof. \quad According to (\ref{5.2}) and the orthogonality of $\{E_\eta
\}$ with respect to (\ref{4.13}) we have
\begin{equation}\label{4.14a}
\langle S_{\eta^*}^{(I,J)}, E_{\eta^*} \rangle_C = ||E_{\eta^*}||_C^2
= {d_{\eta^*}' e_{\eta^+} \over d_{\eta^*} e_{\eta^+}'} 
|| 1 ||_C^2,
\end{equation}
where the second equality is a formula in \cite{BF98b}.
But since the weight function in (\ref{4.13}) is symmetric
$$
\langle S_{\eta^*}^{(I,J)}, E_{\eta^*} \rangle_C =
{1 \over |{\cal O}_{I,J}|} \langle
S_{\eta^*}^{(I,J)}, {\cal O}_{I,J} E_{\eta^*} \rangle_C =
{a_{\eta^*}^{(I,J)} \over |{\cal O}_{I,J}|}
|| S_{\eta^*}^{(I,J)} ||^2_C =
{{\cal M}_{I,\eta^*} \over |{\cal O}_{I,J}|}
{d_{\eta^{*(-,0)}}' d_{\eta^{*(-,0)}} \over
d_{\eta^*}' d_{\eta^{*(-,-)}}' } || S_{\eta^*}^{(I,J)} ||_C^2.
$$
Equating the two results gives (\ref{M.8}).
\hfill $\square$

We are now in a position to present our constant term identity.

\begin{prop}
We have
\begin{eqnarray}\label{CT}
&& {\rm CT} \Big \{ \prod_{i=1}^N (1 - x_i)^a \Big ( 1 - 
{1 \over x_i} \Big )^b \Big (
\prod_s \prod_{j \in \tilde{J}_s}
\Big ( 1 - {1 \over x_j} \Big )^{1 - |\tilde{J}_s|}
\Delta_{\tilde{J}_s}(1/x) \Delta_{\tilde{J}_s}(x) \Big )
\Big ( \Delta(1/x) \Delta(x) \Big )^{1/\alpha}
U_{\eta^*}(x) \Big \} \nonumber \\
&& \quad = {e_{\eta^+}' \over d_{\eta^*}'
[-a-b]_{\delta^+}^{(\alpha)}}
\prod_{\vec{\beta} \in R_{J,+}} \Big (
(\overline{\eta^*}_{\vec{\beta}} - 1)/\alpha \Big )
|| S_{\eta^*}^{(I,J)} ||_C^2
{[-b]_{\eta^+}^{(\alpha)}
 \over
[a+1 + (N-1)/\alpha ]_{\eta^+}^{(\alpha)} }  \nonumber \\
&& \qquad \times
\prod_{i=1}^N {\Gamma(1 + (i-1)/\alpha) \over
\Gamma(1 + a + (i-1)/\alpha)}
{\Gamma(1+a+b+(i-1)/\alpha) \over
\Gamma(1+b+(i-1)/\alpha)}.
\end{eqnarray}
\end{prop}

\noindent Proof. 
\quad From the orthogonality of $\{S_{\eta^*}^{(I,J)}\}$ with respect to
(\ref{4.13}) and the definition (\ref{2.I.1}) of
$U_{\eta^*}^{(I,J)}$, it follows from (\ref{N.0})  that
\begin{eqnarray}\label{U.0}
&& {\rm CT} \Big \{ \prod_{i \notin \tilde{J}}\Big (1 - {1 \over
x_i} \Big )^{r-1} 
\Big (
\prod_s \prod_{j \in \tilde{J}_s}
\Big ( 1 - {1 \over x_j} \Big )^{r - |\tilde{J}_s|}
\Delta_{\tilde{J}_s}(1/x) \Delta_{\tilde{J}_s}(x) \Big )
\Big ( \Delta(1/x) \Delta(x) \Big )^{1/\alpha}
U_{\eta^*}(x) \Big \} \nonumber \\
&& \quad = {\alpha^{|\eta|} \over d_{\eta^*}'}
{[1-r]_{\eta^+}^{(\alpha)} \over
[1-r]_{\delta^+}^{(\alpha)}}
\prod_{\vec{\beta} \in R_{J,+}} \Big (
(\overline{\eta^*}_{\vec{\beta}} - 1)/\alpha \Big )
|| S_{\eta^*}^{(I,J)} ||_C^2
\end{eqnarray}

Substitute $r=a+b+1$ ($a \in \mathbb Z^+$) and make the replacement
$\eta^* \mapsto \eta^* + a^N$, where
$\eta^* + a^N = (\eta_1^*+a, \dots, \eta_N^*+a)$. Since
\begin{equation}\label{U.1}
U_{\eta^*+a^N}(x) = x^a U_{\eta^*}(x)
\end{equation}
we see that
\begin{equation}
|| S_{\eta^*+a^N}^{(I,J)} ||_C^2 = || S_{\eta^*}^{(I,J)} ||_C^2,
\end{equation}
while the definition of $\bar{\eta}_i$ in (\ref{1.7b}) shows
$\overline{\eta^*}_{\vec{\beta}} = \overline{\eta^*}_i - 
\overline{\eta^*}_j$ for
some $i,j$ and thus
\begin{equation} 
\overline{\eta^* + a^N}_{\vec{\beta}} = \overline{\eta^*}_{\vec{\beta}}.
\end{equation}
The definition (\ref{5.3}) of $d_{\eta^*}'$ shows
\begin{eqnarray}\label{U.4}
d_{\eta^* + a^N}' & = & d_{\eta^*}' \prod_{j=1}^a \prod_{i=1}^N
\Big ( N - i + \alpha(j + \eta_i^+) \Big ) \nonumber \\
& = & d_{\eta^*}' \alpha^{aN}
{[a+1+(N-1)/\alpha]_{\eta^+}^{(\alpha)} \over
[1+(N-1)/\alpha]_{\eta^+}^{(\alpha)} }
\prod_{i=1}^N {\Gamma(a+1+ (i-1)/\alpha) \over
\Gamma(1+(i-1)/\alpha))} \nonumber \\
& = & d_{\eta^*}' \alpha^{|\eta|+aN} {[a+1+(N-1)/\alpha]_{\eta^+}^{(\alpha)}
\over e_{\eta^+}'}
\prod_{i=1}^N {\Gamma(a+1+ (i-1)/\alpha) \over
\Gamma(1+(i-1)/\alpha))} 
\end{eqnarray}
where the second equality follows upon use of (\ref{se})
and the third upon use of (\ref{4.14'}). It
also follows from (\ref{se}) that
\begin{equation}\label{U.5}
[-a-b]_{\eta^+ + a^N}^{(\alpha)}  = [-b]_{\eta^+}^{(\alpha)}
(-1)^{aN} \prod_{j=1}^N {\Gamma(1+a+b+(j-1)/\alpha) \over
\Gamma(1+b+(j-1)/\alpha)}.
\end{equation}
The formulas (\ref{U.1}) -- (\ref{U.5}) show that the RHS of
(\ref{U.0}) reduces to the RHS of (\ref{CT}) except that there is
an extra factor $(-1)^{aN}$.

On the LHS, after noting
$$
\prod_{i \notin \tilde{J}_s} \Big (1 - {1 \over x_i} \Big )^{r-1}
\prod_s \prod_{j \in \tilde{J}_s}
\Big ( 1 -
{1 \over x_i} \Big )^{r - | \tilde{J}_s|}
= \prod_{i=1}^N \Big ( 1 -
{1 \over x_i} \Big )^{r - 1}
\prod_s \prod_{j \in \tilde{J}_s}
\Big ( 1 -
{1 \over x_i} \Big )^{1 - | \tilde{J}_s|}
$$
and making use of (\ref{U.1}), we obtain the constant term on the LHS
of (\ref{CT}), except for an extra factor of $(-1)^{aN}$.
\hfill $\square$ 

In the case $I=J= \emptyset$ and thus $\eta^* = \eta$, $S_{\eta^*}^{(I,J)}
= E_\eta$, (\ref{CT}) gives
\begin{eqnarray*}
&&{{\rm CT} \{ \prod_{i=1}^N (1 - x_i)^a (1 - 1/x_i)^b
(\Delta(1/x) \Delta(x) )^{1/\alpha} E_\eta(x) \} \over
{\rm CT} \{ \prod_{i=1}^N (1 - x_i)^a (1 - 1/x_i)^b
(\Delta(1/x) \Delta(x) )^{1/\alpha} \} } \\
&& \quad= { e_{\eta^+}' \over d_\eta'}
{|| E_\eta ||_C^2 \over ||1||_C^2}
{[-b]_{\eta^+}^{(\alpha)} \over
[a+1+(N-1)/\alpha]_{\eta^+}^{(\alpha)} } =
E_\eta(1^N) {[-b]_{\eta^+}^{(\alpha)} \over
[a+1+(N-1)/\alpha]_{\eta^+}^{(\alpha)} }
\end{eqnarray*}
where the second equality follows from (\ref{4.14a}) and (\ref{2.28a}).
This is the constant term identity given in \cite{BF98b}. Summing over
appropriate linear combinations of $\eta : |\eta| = \kappa$ so that
the nonsymmetric Jack polynomial $E_\eta$ becomes the symmetric
Jack polynomial $P_\kappa$ gives the original Macdonald-Kadell-Kaneko
identity.

%\bibliographystyle{plain}
%\bibliography{book}

\begin{thebibliography}{10}

\bibitem{BDF00}
T.H. Baker, C.F. Dunkl, and P.J. Forrester.
\newblock Polynomial eigenfunctions of the Calogero-Sutherland-Moser models
  with exchange terms.
\newblock In J.F. van Diejen and L.~Vinet, editors, {\em
  Calogero-Sutherland-Moser models}, CRM Series in Mathematical Physics, pages
  {37--51}. Springer, New York, 2000.

\bibitem{BF97b}
T.H. Baker and P.J. Forrester.
\newblock The {Calogero-Sutherland} model and polynomials with prescribed
  symmetry.
\newblock {\em Nucl. Phys. B}, 492:682--716, 1997.

\bibitem{BF98b}
T.H. Baker and P.J. Forrester.
\newblock Nonsymmetric Jack polynomials and integral kernels.  
\newblock {\em Duke Math. J.}, 95:1--50, 1998.

\bibitem{Du98}
C.F. Dunkl.
\newblock Orthogonal polynomials of types {$A$} and {$B$} and related Calogero
  models.
\newblock {\em Commun. Math. Phys.}, 197:451--487, 1998.

\bibitem{Ka97g}
K.W.J. Kadell.
\newblock The {Selberg-Jack} symmetric functions.
\newblock {\em Adv. Math.}, 130:33--102, 1997.

\bibitem{Ka93}
J.~Kaneko.
\newblock Selberg integrals and hypergeometric functions associated with {Jack}
  polynomials.
\newblock {\em SIAM J. Math Anal.}, 24:1086--1110, 1993.

\bibitem{KY98}
Y.~Kato and T.~Yamamoto.
\newblock Jack polynomials with prescribed symmetry and hole propogator of spin
  Calogero-Sutherland model.
\newblock {\em J. Phys. A}, 31:9171--9184, 1998.

\bibitem{KS97}
F.~Knop and S.~Sahi.
\newblock A recursion and combinatorial formula for Jack polynomials.
\newblock {\em Inv. Math.}, 128:9--22, 1997.

\bibitem{Ma87}
I.G. Macdonald.
\newblock Commuting differential operators and zonal spherical functions.
\newblock In A.M.~Cohen et~al., editor, {\em Algebraic Groups, Utrecht 1986},
  volume 1271 of {\em Lecture Notes in Math.}, pages 189--200. Springer-Verlag,
  Heidelberg, 1987.

\bibitem{NW00}
A.~Nishino and M.~Wadati.
\newblock Bosonic and fermionic eigenstates for generalized {S}utherland
  models.
\newblock {\em J. Phys. A}, 33:3795--3807, 2000.

\bibitem{Op98}
E.M. Opdam.
\newblock Lectures on Dunkl operators.
\newblock math.RT/9812007 

\bibitem{Op95}
E.M. Opdam.
\newblock Harmonic analysis for certain representations of graded Hecke
  algebras.
\newblock {\em Acta Math.}, 175:75--121, 1995.

\bibitem{Sa96}
S.~Sahi.
\newblock A new scalar product for nonsymmetric Jack polynomials.
\newblock {\em Int. Math. Res. Not.}, 20:997--1004, 1996.

\end{thebibliography}

\section*{Acknowledgements}
This work was supported by the Australian Research Council. Part of
this work was written up while PJF was a participant in 
the Newton Institute 2001 semester on Symmetric functions
and Macdonald polynomials.

\end{document}